\documentclass[leqno,12pt]{article}
\usepackage{etoolbox} \usepackage{amsfonts}
\usepackage[margin=1in]{geometry}
\usepackage{times}
\usepackage{amsthm,amssymb,stmaryrd}
\usepackage{amsmath,xypic}

\usepackage[usenames,dvipsnames]{color}
\usepackage[authoryear,sort&compress]{natbib}

\usepackage[colorlinks,citecolor=blue,colorlinks=true,hyperindex, citecolor=cyan, urlcolor=cyan]{hyperref}

\usepackage{epigraph}
\usepackage{enumitem}
\usepackage{mathrsfs}
\DeclareMathAlphabet{\mathscrbf}{OMS}{mdugm}{b}{n} \newcommand{\be}{\begin{equation}}
\newcommand{\ee}{\end{equation}}
\newcommand{\bes}{\begin{equation*}}
\newcommand{\ees}{\end{equation*}}
\newcommand{\bea}{\begin{eqnarray}}
\newcommand{\eea}{\end{eqnarray}}
\newcommand{\beas}{\begin{eqnarray}}
\newcommand{\eeas}{\end{eqnarray}}
\newcommand{\ben}{\begin{note}}
\newcommand{\een}{\end{note}}
\newcommand{\bexl}{\vskip0.1em\noindent\hrulefill\vskip1em\begin{ExerciseList}}
\newcommand{\eexl}{\end{ExerciseList}\hrulefill}

\newcommand{\bthm}{\begin{theorem}}
\newcommand{\ethm}{\end{theorem}}
\newcommand{\bpro}{\begin{prop}}
\newcommand{\epro}{\end{prop}}
\newcommand{\bcor}{\begin{corollary}}
\newcommand{\ecor}{\end{corollary}}
\newcommand{\bcon}{\begin{conjecture}}
\newcommand{\econ}{\end{conjecture}}
\newcommand{\bp}{\begin{proof}}
\newcommand{\ep}{\end{proof}}
\newcommand{\blem}{\begin{lemma}}
\newcommand{\elem}{\end{lemma}}
\newcommand{\bn}{\begin{note}}
\newcommand{\en}{\end{note}}
\newcommand{\benum}{\begin{enumerate}}
\newcommand{\eenum}{\end{enumerate}}
\newcommand{\bed}{\begin{defn}}
\newcommand{\eed}{\end{defn}}
\newcommand{\brem}{\begin{remark}}
\newcommand{\erem}{\end{remark}}

\newcommand{\btik}{\begin{tikzpicture}\begin{axis}[scale=0.5,axis y line=center, axis x line=middle]}
\newcommand{\etik}{\end{axis}\end{tikzpicture}}

\let\into=\hookrightarrow
\let\mapsto=\longmapsto

\newcommand{\upperRomannumeral}[1]{\uppercase\expandafter{\romannumeral#1}}

 \usepackage{tikz-cd}
\usepackage{stmaryrd}

\newtheorem{theorem}[equation]{Theorem}      \newtheorem{lemma}[equation]{Lemma}          \newtheorem{corollary}[equation]{Corollary}  \newtheorem{proposition}[equation]{Proposition}

\theoremstyle{definition}

\theoremstyle{definition}
\newtheorem{defn}[equation]{Definition}
\theoremstyle{remark}

\theoremstyle{definition}
\newtheorem{remark}[equation]{Remark}

\numberwithin{equation}{section}

\let\cite=\citep

\let\into=\hookrightarrow
\let\isom=\simeq

\let\tensor=\otimes

\newcommand{\bF}{{\bar{F}}}

\newcommand{\C}{{\mathbb C}}

\newcommand{\End}{\rm{End}}

\newcommand{\gal}{{\rm Gal}}

\newcommand{\mydot}{{\small{\bullet}}}
\newcommand{\N}{\mathcal{N}}

\newcommand{\Q}{{\mathbb Q}}

\newcommand{\spec}{{\rm Spec}}
\newcommand{\Spec}{{\rm Spec}}

\newcommand{\Z}{{\mathbb Z}}

\renewcommand{\int}{\operatorname{int}}
\renewcommand{\O}{{\mathcal O}}

\newcommand{\mapright}[1]{{\xymatrix{{}\ar[r]^{#1}&{}}}}

\renewcommand{\bpro}{\begin{proposition}}
	\renewcommand{\epro}{\end{proposition}}

\let\mathcal=\mathscr

\title{Algebraization of Mochizuki's anabelian variation of ring structures, perfectoid geometry and formal groups}
\author{Kirti Joshi}

\begin{document}

\maketitle

\epigraphwidth0.65\textwidth
\epigraph{Ry\=okan! how nice to be like  a fool\\
	$\qquad$ for then one's Way is grand beyond all measure}{(Master) Tainin Kokusen (to Ry\=okan Taigu) \cite{ryokan}}

\begin{abstract}
Let $M$ be a multiplicative monoid with identity. Then I show that there is a universal one dimensional formal group law equipped with an action of $M$. If $M$ is $p$-perfect (i.e. $m\mapsto m^p$ is an isomorphism for some prime number $p$) then the universal $M$-formal group law comes equipped with a natural Frobenius endomorphism.  There are a number of concrete applications of this result. If $K$ is a $p$-adic field and $\O=\O_K$ is the multiplicative monoid of the ring of integers of $K$, then there is a universal formal group (over a suitable (non-zero) ring) which is equipped with an action of the multiplicative monoid $\O$. Lubin-Tate formal groups arise from this universal monoid formal group law. This has applications to Mochizuki's anabelian ideas: if two p-adic fields have isomorphic absolute Galois groups then they have isomorphic multiplicative monoids $\O$ (but possibly non-isomorphic ring structures). The existence of the universal monoid formal group law for the monoid $\O$ implies that  the additive structures of a ring  can be interpolated  into a universal algebraic family (while keeping the multiplicative structure of the ring fixed). Here is another important example covered by my result:  let $R$ be a perfectoid ring and let $R^\flat$ be its tilt and the multiplicative monoid  $R^\flat$ of $R^\flat$. Then there exists a universal monoid formal group law for this monoid which interpolates the additive structures of untilts with tilt $R^\flat$. Thus in some sense one has a unified approach to various phenomenon which are well-known in anabelian geometry and in perfectoid geometry. These results also provide a natural number field version of Fontaine's fundamental ring $A_{inf}$ of $p$-adic Hodge Theory (Subsection~\ref{ss:global-ainf}).
\end{abstract}

\tableofcontents 
\newcommand{\act}{\curvearrowright}
\newcommand{\lmp}{{\Pi\act\Ot}}
\newcommand{\lmpi}{{\lmp}_{\int}}
\newcommand{\lmpf}{\lmp_F}
\newcommand{\Om}{\O^{\times\mu}}
\newcommand{\Omf}{\O^{\times\mu}_{\bF}}
\renewcommand{\N}{\mathbb{N}}
\newcommand{\yoga}{Yoga}

\newcommand{\sM}{\mathcal{M}}

\newcommand{\sF}{\mathscr{F}}
\newcommand{\stm}{M^\circledast}
\newcommand{\cdas}{\circledast}
\newcommand{\okm}{\O_K^\triangleright}
\newcommand{\tr}{\triangleright}
\newcommand{\pow}[2]{#1\llbracket#2\rrbracket}
\newcommand{\isomar}{\overset{\sim}{\longrightarrow}}
\newcommand{\benumlab}{\benum[label={\bf(\arabic{*})}]}
\newcommand{\iut}{\cite{mochizuki-iut1,mochizuki-iut2,mochizuki-iut3,mochizuki-iut4}}
\newcommand{\br}{\begin{remark}}
\newcommand{\er}{\end{remark}}

\newcommand{\ocs}{\O^\circledast}
\newcommand{\osm}{\O^*}
\newcommand{\bosm}{\bar{\O}^*}
\newcommand{\os}{\osm}
\newcommand{\bos}{\bar{\O}^*}
\newcommand{\ot}{\O^\tr}
\renewcommand{\bot}{\bar{\O}^\tr}
\newcommand{\omu}{\O^\mu}
\newcommand{\bomu}{\bar{\O}^\mu}
\newcommand{\sFm}{\sF_{\okm}}
\newcommand{\sFs}{\sF_{\O^*}}
\newcommand{\sFmu}{\sF_{\omu}}
\newcommand{\sFk}{\sF_{\ks}}
\newcommand{\lks}{L_{\ks}}
\newcommand{\los}{L_{\osm}}
\newcommand{\lomu}{L_{\omu}}
\newcommand{\lot}{L_{\O^\tr}}
\newcommand{\ks}{\mathscr{K}^*}
\newcommand{\bks}{\bar{\ks}}
\newcommand{\bK}{\bar{K}}

\newcommand{\kxt}{\mathscr{K}^*}
\newcommand{\kx}{\mathscr{K}_X}
\newcommand{\ky}{\mathscr{K}_Y}
\newcommand{\kss}{\mathscr{K}^\circledast}
\newcommand{\kxg}{\mathscr{K}^*_{\bf g}}
\newcommand{\kxs}{\mathscr{K}^\circledast_X}
\newcommand{\kgs}{\mathscr{K}^\circledast_{\bf g}}
\newcommand{\gah}[1]{\hat{\mathbb{G}}_{a/{#1}}}

\newcommand{\bbW}{\mathbb{W}}
\newcommand{\bbWl}{\bbW_L}
\newcommand{\ainf}{A_{inf}}
\newcommand{\ainfl}{\mathbb{A}_{inf,L}}
\newcommand{\bL}{\bar{L}}
\newcommand{\cpt}{\C_p^\flat}
\newcommand{\invlim}{\varprojlim}

\newcommand{\topics}{\cite{mochizuki-topics1,mochizuki-topics2,mochizuki-topics3}}

\numberwithin{equation}{subsection}

\section{Introduction}
\subsection{The key question}
In this paper, the term  \emph{$p$-adic field} will mean  a finite extension of $\Q_p$ for some prime number $p$. Let $G_K$ be the absolute Galois group of $K$ for some choice of an algebraic closure of $K$. Let $G_K^\mydot\subset G_K$ be the upper numbering ramification filtration of $G_K$ (\cite{serre1979-local-fields}). In \cite{mochizuki-local-gro}, it is shown that a $p$-adic field $K$ can be recovered from the pair $(G_K,G_K^\mydot)$ consisting of the topological $G_K$ equipped with its filtration by inertia subgroups. On the other hand, from  \cite{yamagata76}, \cite{jarden79} one knows that, for every prime number $p$, there exists non-isomorphic $p$-adic fields $K_1$ and $K_2$   such that one has an isomorphism topological groups:
\benumlab
\item $G_{K_1}\isom G_{K_2}$,
\item $K^*_1\isom K^*_2$.
\eenum
So one could say that, for a $p$-adic field $K$, the multiplicative monoid $K^*\cup \{ 0\}$ can be equipped with many non-isomorphic ring structures i.e. there exists many different ring structures on this set in which the multiplicative structure remains fixed (up to an isomorphism) or one can loosely say that the additive structure deforms while the multiplicative structure remains fixed. This phenomenon arose in classical works and more recent works  in anabelian geometry, especially \topics, \iut, and so it is natural to refer to it as the anabelian fluidity of the ring structures. This fluidity is an integral part of the theory developed by Mochizuki in \cite{mochizuki-theta,mochizuki-topics1, mochizuki-topics2,mochizuki-topics3,mochizuki-iut1,mochizuki-iut2,mochizuki-iut3,mochizuki-iut4}.

However, the predominant view (and toolkit) in the aforementioned works is that of group theory i.e. group theory surrounding the absolute galois groups or \'etale fundamental groups. This group-theoretic view is extremely unsatisfactory from an algebro-geometric point of view. The question which occurred to me is this:

\begin{center}
\textit{Can this (anabelian) fluidty of ring structures be algebraized in some way?}
\end{center}

\subsection{The key insight}
I began  thinking of this question in Kyoto (Spring 2018). The idea, which I elaborate here, occurred to me in a lecture by Michael Hopkins at the Arizona Winter School (2019). In one of his lectures, Hopkins narrated an anecdote about Daniel Quillen's discovery (\cite{quillen69}) of the role of formal groups in topological cohomology theories: in particular Quillen's assertion (to Hopkins) that ``as addition rule for Chern classes fails to hold, it must therefore fail in the worst possible way--namely by means of a formal group'' (I am paraphrasing both Hopkins and Quillen here). 

It was immediately clear to me after listening to Hopkin's anecdote, that the problem of fluidity of additive structures of a $p$-adic field is also a problem of deforming the tautological addition law 
\be\label{eq:tauto}(x,y)\mapsto x+y\ee of a field and that a formal group law
$$(x,y)\mapsto x+y+ \text{higher order terms}$$
provides the required deformation of the tautological addition law \eqref{eq:tauto}. 
\begin{center}
	\textit{This algebraization is the central discovery of this paper.} 
\end{center}
Importantly,  not only can this problem be algebraized, my approach also provides  a unified view of both anabelian fluidity of ring structures and similar fluidity of ring structures (with fixed multiplicative structures) encountered in perfectoid geometry \cite{scholze12-perfectoid-ihes},  \cite{kedlaya-liu15} and \cite{fargues-fontaine}, and also provides similar objects in the global arithmetic case of varieties over number fields  and also provides the global ring $\ainfl$ which is the natural analog of Fontaine's fundamental ring $\ainf$ of $p$-adic Hodge Theory (\cite{fontaine94a}).

This paper is a record of these ruminations. 

\subsection{Results  proved in this paper}
I will assume that readers are familiar with basic aspects of the theory of formal groups as documented in  \cite{hazewinkel-book} and I will freely use notation of that book. The new idea introduced here is that of a monoid formal group law (Definition~\ref{de:monoid-frml-grp}). One dimensional Lubin-Tate formal groups provide non-trivial examples of monoid formal group laws. In Theorem~\ref{th:universal-M-formal-group}, I show that for any given commutative monoid $M$ (written multiplicatively), there exists a universal (monoid) $M$-formal group law. If the monoid $M$ is $p$-perfect (for some prime number $p$) (this means $m\mapsto m^p$ is an isomorphism of $M$) then one obtains (Theorem~\ref{th:frob-exist}) a Frobenius morphism on the universal object provided by Theorem~\ref{th:universal-M-formal-group}. 

Let $K$ be a $p$-adic field and let $\O_K^\tr$ be the multiplicative monoid of non-zero elements of the ring of integers $\O_K$ of $K$. In Theorem~\ref{th:converse-all-ring-structures}, I prove that every commutative ring $R$ of characteristic zero (with unity) such that $R^\tr\isom \O^\tr_K$ arises from the universal $\okm$-formal group in a manner prescribed by Proposition~\ref{pr:formal-groups-and-ring-struct} and more importantly, that every discrete valuation ring $R$ of characteristic zero with a finite residue field such that $R^\tr\isom \O^\tr_K$ arises from a non-trivial (i.e. not isomorphic to the additive group) $\okm$-formal group and hence in particular from the universal $\okm$-formal group in the manner prescribed by Proposition~\ref{pr:formal-groups-and-ring-struct}. \emph{So  one could say that the universal $\okm$-formal group gives rise to all ring structures of interest on $\okm\cup\{ 0\}$}.

In Section~\ref{se:geometric-applications}, I provide two important geometric applications. In Subsection~\ref{ss:geometric-case}, Theorem~\ref{th:universal-M-formal-group}, Theorem~\ref{th:frob-exist} are applied to construct deformations of monoids provided by a geometrically connected smooth quasi-projective variety in Theorem~\ref{th:universal-oxt-formal-group}. In particular, one also obtains deformations of the sort considered in \cite{mochizuki-iut2,mochizuki-iut3}. 

In Subsection~\ref{ss:perfectoid-case}, Theorem~\ref{th:universal-M-formal-group}, Theorem~\ref{th:frob-exist} are applied to perfectoid geometry. In  Theorem~\ref{th:perfectoid-case}, I demonstrate how perfectoid rings, fields and more generally, perfectoid Huber pairs (\cite{scholze12-perfectoid-ihes}, \cite{kedlaya-liu15}, \cite{fargues-fontaine}) also relate to the idea of fluidity of ring structures--especially, untilts with a fixed tilt are bound together by a universal addition law. The special case of  perfectoid fields and the relationship with \cite{fargues-fontaine} is detailed in \cite{joshi-teich-def}. 

In subsection~\ref{ss:global-ainf}, using the results of Subsection~\ref{ss:geometric-case} and Subsection~\ref{ss:perfectoid-case},  I sketch the construction of the global (i.e. number field) analog $\ainfl$ of the  fundamental ring of $p$-adic Hodge Theory $\ainf$ (introduced by Jean-Marc~Fontaine in \cite{fontaine94a}).

Hence this work not only provides a natural and unified view of two apparently disparate themes in modern algebraic geometry, but also demonstrates that anabelian geometry and perfectoid geometry are inseparably linked together and much insight may be gained by understanding this union. Some aspects of this are already reflected in my work on Arithmetic Teichmuller Spaces. 

\subsection{Acknowledgments} I began thinking about the problem addressed in this paper while I was on a sabbatical (Spring 2018) and enjoying the support and hospitality of Research Institute of Mathematical Sciences (RIMS), Kyoto and I thank RIMS for the same. I thank Shinichi Mochizuki for conversations around many topics of common interest during my visit. I thank Yuichiro Hoshi for conversations about Anabelian Reconstruction Theory. I thank Machiel van Frankenhuijsen for many conversations related to Mochizuki's work and the $abc$-conjecture. I thank Taylor Dupuy and Anton Hilado for some conversations surrounding Mochizuki's works and for pointing out some typos in an earlier version of this paper.  Finally it is a pleasure to thank Michael Hopkins for his insightful lectures at the Arizona Winter School 2019 and  the organizers of the School for an excellently organized meeting.  

\newcommand{\mybrackets}[1]{\llbracket#1\rrbracket}

\section{Formal groups with monoid actions}\label{se:formal-groups}
In this paper a \textit{$p$-adic field} will be a finite extension of $\Q_p$ for some prime number $p$. A \textit{perfectoid field} will be as defined in \cite[Definition 3.1]{scholze12-perfectoid-ihes}. If $K$ is a field, I will write $G_K$ for its absolute Galois group for some choice of a separable closure of $K$. If $K_1,K_2$ are fields then I will say, following \cite[Definition 2.1.1]{joshi-anabelomorphy}, that $K_1,K_2$ are \textit{anabelomorphic fields} if one has an isomorphim $G_{K_1}\isom G_{K_2}$ of topological groups. In most of this  paper, $K$ will be a $p$-adic field. Later I will also allow $K$ to be a perfectoid field. For more on anabelomorphisms of $p$-adic fields, perfectoid fields or number fields, the reader is referred to \cite{joshi-anabelomorphy}. 
\subsection{Monoid formal groups}
Let $R$ be a commutative ring with a unit. Write $R^*$ for the group of units in $R$, and if $R$ is an integral domain, let $R^\tr=R-\{0\}$ be the multiplicative monoid of non-zero elements in $R$. If $R$ is any ring let $R^\cdas$ be the multiplicative monoid of $R$.
I will work with formal groups over $R$. All formal groups I deal with will be one dimensional commutative formal groups over the chosen ring $R$ (I will not repeat these assumptions again). Let $\sF$ be such a formal group. So $\sF(x,y)\in R\mybrackets{x,y}$ is a formal power series in $x,y$ with coefficients in $R$ with certain properties (see \cite[Chap. 1, 1.1.1--1.1.3]{hazewinkel-book}). 

\begin{defn}\label{de:monoid-frml-grp}
Let $M$ be a commutative monoid with identity (\emph{always written multiplicatively in this paper}). I say that $\sF$ is a \emph{formal group with an action of the monoid $M$} (or simply say that \emph{$\sf$ is an $M$-formal group over $R$}) if one has a homomorphism of multiplicative monoids $$\iota_M:M\to \End_R(\sF).$$ 
\end{defn}
For $m\in M$, write $[m]\in \End_R(\sF)$ for the endomorphism $\iota_M(m)$ corresponding to $m\in M$. Note that this means $[m]$ is a formal series over $R$ in some auxiliary variable, say $x$, and when one wants to emphasize this, I will write $[m](x)$ for an auxiliary variable $x$. So $[m]$ is of the form 
\be [m](x)=\alpha_1(m)\cdot x+\sum_{i\geq 2} \alpha_i(m)x^i.\ee

\blem
The mapping  given by $m\mapsto \alpha_1(m)$ is a homomorphism of multiplicative monoids 
$$ \alpha_1: M\to R.$$
\elem
\bp 
This is easily verified using the fact that $m\mapsto[m]$ is a homomorphism of multiplicative monoids.
\ep

\brem 
The  notion $M$-formal groups arose from a natural question which occurred to me while reading Mochizuki's ideas about anabelian changes of ring structures in \iut.  Let  $K_1$ and $K_2$ be two anabelomorphic $p$-adic fields i.e. $K_1$ and $K_2$ are $p$-adic fields with  topologically isomorphic absolute Galois groups. Then by well-known classical results (see \cite[Theorem 2.4.3]{joshi-anabelomorphy} for instance), one has an isomorphism of the topological multiplicative monoids $K_1^*\isom K_2^*$. On the other hand these $p$-adic fields need not be isomorphic. So one has two possibly non-isomorphic fields whose respective multiplicative monoids are isomorphic. The predominant view in anabelian geometry is that of group theory i.e. group theory surrounding the absolute galois groups or \'etale fundamental groups. This group-theoretic view is extremely unsatisfactory from an algebro-geometric point of view. The question which occurred to me is whether this phenomenon could be algebraized in some way. This algebraization is the discovery of this paper. 
\erem

\subsection{Strict $M$-formal groups}
Let $M$ be a monoid in the convention of Definition~\ref{de:monoid-frml-grp}. Let $$\stm=M\cup \{ 0\}$$  with multiplication $0\cdot m=0=m\cdot 0$ for all $m\in M$. Note that $\stm$ is obviously a commutative monoid with this rule.  Suppose $\sF$ is an $M$-formal group over $R$. Then clearly $\iota:M\to \End_R(\sF)$ extends to a homomorphism of multiplicative monoids $$\iota^\cdas:\stm \to \End_R(\sF) \text{ with } \iota^\cdas(0)=0\in\End_R(\sF).$$

\bpro\label{pr:formal-groups-and-ring-struct}
Suppose $\sF$ is an $M$-formal group over $R$. Suppose that $\iota_M$ induces an isomorphism of multiplicative monoids $\iota^\cdas_M:\stm\isomar \End_R(\sF)$.
Then there exists a unique structure of a commutative ring on $\stm$ given by $\sF$, such that $M$ is the multiplicative monoid of non-zero elements.   The additive structure of $\stm$ is given by
\begin{align}
[m_1+m_2](x) & =  \sF([m_1](x),[m_2](x)) \\
[m_1\cdot m_2]  & =  [m_1]\cdot [m_2].
\end{align}
\epro

\bp 
The proof is as follows: $\End_R(\sF)$ is a ring, and if $\iota^\cdas:\stm\to \End_R(\sF)$ is a bijection of multiplicative monoids, then $\stm$ carries a ring structure induced from the target of $\iota^\cdas$ and in this ring structure $M$ is multiplicative monoid of non-zero elements. As $\stm$ is a commutative monoid and as $\iota^\cdas_M:\stm\isomar \End_R(\sF)$ is an isomorphism of multiplicative monoids, this ring structure  on $\stm$ is that of a commutative ring.  This ring structure is given by the stated formulae which simply state the rules which hold in the target of $\iota^\cdas$. So for example, to add $m_1,m_2\in M$, one adds the endomorphisms $[m_1],[m_2]$ of our formal group $\sF$ and the result gives $m_1+m_2$ (uniquely). 
\ep

This proposition motivates the following definition. 

\begin{defn} Let $\sF$ be an $M$-formal group law over a ring $R$. I say that $\sF$ is a \emph{strict $M$-formal group over $R$} if $\iota_M^\cdas$ induces an isomorphism of multiplicative monoids $\iota^\cdas_M:\stm\isomar \End_{R}(\sF)$. 
\end{defn}

\brem 
As one sees from Lemma~\ref{le:lubin-tate-strict}, Lubin-Tate formal group laws are examples of monoid formal group laws which are strict.
\erem 

\bpro 
If $\sF$ is any strict $M$-formal group law over a ring $R$, then $\End_{R}(\sF)$ is a commutative ring.
\epro
\bp 
Since $\sF$ is a strict $M$-formal group law, the homomorphism $M^\cdas \to \End_{R}(\sF)$ is an isomorphism of monoids. Since $M$ is a commutative monoid, so is $M^\cdas$ and hence so is $\End_{R}(\sF)$. Thus the multiplicative monoid $\End_{R}(\sF)$ is commutative. Thus the ring $\End_{R}(\sF)$ is a commutative ring.
\ep

\subsection{Existence of formal groups with monoid actions}
Let $\Q[M]$ be the monoid ring of a monoid $M$ over $\Q$.
Let $R\supset\Q[M]$ be a commutative ring. Then I claim that there always exists an $M$-formal group over $R$.
\bthm\label{th:exist-M-formal} Let $M$ be a commutative monoid and let $\Q[M]$ be the monoid ring over $\Q$.
Let $R\supseteq\Q[M]$ be a commutative ring.  Then there exists an $M$-formal group $\sF$ over $R$.
\ethm
\bp
One uses the fact that over any $\Q$-algebra, any one-dimensional commutative formal group is isomorphic to the additive formal group with the group law $$\sF(x,y)=x+y.$$ Then it is clear that $$[m](x)=m\cdot x$$ defines an endomorphism of our additive group as our ring $R\supset \Q[M]\supset M$ and it is clear that $m\mapsto [m]\in \End_R(\sF)$  provides  a homomorphism of multiplicative monoids. Thus an $M$-formal group exists over any $R\supset\Q[M]$.
\ep

Here is another constructive proof of Theorem~\ref{th:exist-M-formal}.

\bp[Second proof of Theorem~\ref{th:exist-M-formal}]This proof closely follows the proof of existence of formal groups (see \cite{hazewinkel-book}).
Let $f(T)\in \pow{R}{T}$ be a power series in $T$ with $f(T)=T\bmod{T^2}$. Let $g(T)\in \pow R T$ be the unique power series such that $f(g(T))=T=g(f(T))$. Such a $g$ always exists as $\Q\subseteq\Q[M]\subseteq R$.

Now define $$\sF(x,y)= g(f(x)+f(y)),$$
and for any $m\in M$, let $[m]$ be the endomorphism of $\sF$ defined by
$$
[m](x):= g(m\cdot f(x)).
$$
Note that $m\cdot f(x)$ makes sense as $m\in R\supset\Q[M]$ and $f(x)\in \pow R x$.
\ep

\brem  In Theorem~\ref{th:universal-M-formal-group} I show that there exists a universal $M$-formal group which gives rise to all the $M$-formal group laws over any commutative ring. 
\erem

\subsection{Existence of strict $\O^\tr_K$-formal groups and rings structures on $\O^\tr_K\cup\{0\}$}
Now let me explain why the notion $M$-formal groups is germane  to Mochizuki's ideas in \topics, \iut. Let us take  our monoid $M$ to be $M=\okm$ the multiplicative monoid of non-zero elements of the ring of integers of a $p$-adic field $K$. Let $K_1$ and $K_2$ be two anabelomorphic $p$-adic fields i.e. $K_1$ and $K_2$ have topologically isomorphic absolute Galois groups. Then by \cite[Theorem 2.4.3]{joshi-anabelomorphy}, the multiplicative monoids $\O_{K_1}^*\isom \O_{K_2}^*$ i.e. the multiplicative monoid $\O_K^*$ is an amphoric algebraic structure associated to a $p$-adic field.

Let us record the following which will be used in the sequel:

\bpro\label{le:lubin-tate-strict}
Let $K$ be any $p$-adic field. Then there exists a strict $\okm$-formal group law over $\O_K$.  
\epro

\bp 
Indeed by the central result of Lubin-Tate Theory there exists  a one-dimensional Lubin-Tate formal group $\sF$ over $\O_K$ such that one has an isomorphism of rings  (see \cite[Section 20.1.21]{hazewinkel-book}) $$\O_K\isom\End_{\O_K}(\sF).$$ Then the restriction of this homomorphism of rings to the multiplicative monoid $\okm$ makes $\sF$ into a strict $\okm$-formal group over $\O_K$. Hence there exist natural, strict $\okm$-formal group laws over $\O_K$.
\ep

Proposition~\ref{pr:formal-groups-and-ring-struct} and Lemma~\ref{le:lubin-tate-strict} can be applied to our context and let me record the preceding discussion in the following: 

\bthm\label{th:strict-formal-groups-for fields} 
Let $G$ be the absolute Galois group of some $p$-adic field $K$. Let $\O^\tr\isom\okm$ be the multiplicative monoid constructed from $G$. \benumlab
\item Any strict $\okm$-formal group law $\sF$ over $R$ equips the multiplicative monoid $\O^\cdas=\O^\tr\cup\{ 0\}$ with the structure of a commutative ring with $\O^\tr$ as the multiplicative monoid of non-zero elements.
\item Suppose $K_1,K_2$ are anabelomorphic $p$-adic fields. 
\benum 
\item Then through the isomorphism $\O^\tr_{K_1}\isom \O^\tr_{K_2}$, any  $\O^\tr_{K_1}$-formal group over any ring $R$ can be viewed as an $\O^\tr_{K_2}$-formal group over $R$ and vice versa.
\item Any strict $\O^\tr_{K_1}$-formal group over $\O_{K_2}$ equips $\O^\cdas_{K_1}$ with a ring structure (in general distinct from its given ring structure) provided by the given formal group law.
\item More precisely any strict $\O_{K_1}^\tr$-formal Lubin-Tate group  over $\O_{K_2}$ equips $\O^\cdas_{K_1}$ with a ring structure isomorphic to $\O_{K_2}$. 
\eenum
\eenum
\ethm

\section{Existence of Universal $M$-Formal Groups}
\subsection{The existence theorem}
The above theorem shows very explicitly that the fluidity of the ring structures on $\O_K$ arises from the fluidity of strict $\O_K^\tr$-formal group laws and in this sense fluidity  of  additive structures of $\O_K^\tr\cup \{0 \}$ is encoded in the fluidity of certain formal group laws.

This leads one to suspect that fluidity of local additive structures might be a reflection of existence of universal $\O^\tr$-formal group law over some suitably universal ring. This hope is realized  in the theorem proved below. 

\bthm\label{th:universal-M-formal-group}
Let $M$ be any commutative monoid. Then there exists a universal $M$-formal group over a ring $L_M$ and homomorphism of multiplicative monoids $M\to L_M$ with the following universal property. If $\tilde{\sF}$ is an $M$-formal group over some ring $R$ then there is a ring homomorphism $f:L_M\to R$ and $\tilde{\sF}(x,y)$ is obtained by applying $f$ to the coefficients of $\sF(x,y)$. Conversely any homomorphism of rings $L_M\to R$ provides an $M$-formal group on $R$. 
\ethm

\bp The idea of the proof is similar to the proof of the existence of the universal formal group law in \cite[Chapter 1]{hazewinkel-book}. Let $L_M$ be the ring which is being sought. The assertion will be proved if one can prove the existence of a power series in two variables  $$\sF(x,y)=x+y+\sum_{i+j\geq 2} c_{i,j}x^iy^j\in L_M\mybrackets{x,y}$$
with coefficients in this ring $L_M$.  This power series $\sF(x,y)$ is required to satisfy the following list of properties which make it into a formal group over $L_M$:

\benumlab
\item $\sF(x,y)=\sF(y,x)$,
\item $\sF(x,\sF(y,z))=\sF(\sF(x,y),z)$,
\eenum
and the following list of properties which make $\sF$ into an $M$-formal group over $L_M$:
\benum[label={\bf(3)}]
\item for every $m\in M$ an endomorphism $g_m(x)=m\cdot x+\sum_{i\geq 2}d_{m,i}x^i$ such that
\benum
\item $g_1(x)=x$
\item $g_m(\sF(x,y))=\sF(g_m(x),g_m(y)),$
\item $g_{m_1}(g_{m_2}(x))=g_{m_2}(g_{m_1}(x))$
\eenum
\eenum
Let $\Z[M]$ be the monoid ring of $M$ over $\Z$. The idea of the proof is to start with an arbitrary formal power series 
$$\sF(x,y)=x+y+\sum_{i+j\geq 2} c_{i,j}x^iy^j$$
and for each $m\in M$, a formal power series
$$g_m(x)=m\cdot x+\sum_{i\geq 2}d_{m,i}x^i$$
with $g_1(x)=x$ and all these 
with arbitrary coefficients. Now let us write down all the relations which must hold for $\sF$ and $\{g_m(x)\}_{m\in M}$ to be an $M$-formal group and then declare all the coefficients to be variables and eliminate all the obstructions to the relations.
The first property is commutativity
which says $$\sF(x,y)-\sF(y,x)=\sum_{i+j\geq 2}(c_{i,j}-c_{j,i})x^iy^j$$ this gives the condition on the coefficients
$$c_{i,j}=c_{j,i}$$
for commutativity to hold.
Next associativity gives 
$$\sF(x,\sF(y,z))-\sF(\sF(x,y),z)=\sum_{i,j,k}P_{i,j,k}x^{i}y^jz^k,$$ where $P_{i,j,k}$ are polynomials in the coefficients of $\sF$ and hence, in particular, associativity holds if and only if the polynomials $P_{i,j,k}=0$.
Next let us understand the condition that for every $m\in M$ we have an endomorphism of our formal group. This condition means the following: for every $m\in M$, one wants a power series $g_m(x)=m\cdot x+\sum_{i\geq 2}d_{m,i}x^i$ which satisfies
$$g_m(\sF(x,y))-\sF(g_m(x),g_m(y))=\sum_{i}Q_{m,i}x^i$$
where $Q_{m,i}$ are polynomials in the coefficients of the power series $g_m,\sF$.
For instance for $1\in M$ (our monoids are written multiplicatively) one takes $g_1(x)=x$.
These power series $\{g_m\}_{m\in M}$ are required to satisfy the property of composition of endomorphisms:
$$g_{m'}(g_m(\sF(x,y))-g_{m\cdot m'}(\sF(x,y)))=\sum_{i}Z_{m,m',i}x^i,$$
where $Z_{m,m',i}$ is a polynomial in coefficients of $g_m,g_{m'},\sF$
and are required to satisfy commutativity of endomorphisms:
$$g_{m'\cdot m}(x)-g_{m\cdot m'}(x)=\sum_{i\geq 2} (d_{m\cdot m',i}-d_{m'\cdot m,i})x^i.$$
So now let $$L_M=\Z[M][\{c_{i,j}\},\{d_{m,i}\}]/I$$ 
where $I$ is the ideal 
$$I=\left(c_{i,j}-c_{j,i},P_{i,j,k},Q_{m,i},Z_{m,m',i},d_{m\cdot m',i}-d_{m'\cdot m,i}\right).$$
in which all the indices run over all the relevant indexing sets.
Then $\sF$ is clearly a formal group over $L_M$ with a family of endomorphisms $m\mapsto g_m$ for all $m\in M$. Since $g_m(x)=g_{m'}(x)$ if and only if $d_{m,i}=d_{m',i}$ for all $i\geq 2$, so $m\mapsto g_m$  is an injective homomorphism of monoids making $\sF$ into an $M$-formal group over $L$  with $\iota_M:M\to\End_L(\sF)$ given by $m\mapsto g_m$.

Now let me show that this formal group law has the following universal property: given an $M$-formal group  $\tilde{\sF}$ over a ring $R$, there exists a homomorphism $f:L_M\to R$ which provides $\tilde{\sF}$ by applying $f$ to the coefficients of $\sF$. 

To prove this, suppose that $\tilde{\sF}$ is an $M$-formal group over $R$ given by $$\tilde{\sF}(x,y)=x+y+\sum_{i,j} \tilde{c}_{i,j}x^iy^j$$ with $\tilde{c}_{i,j}$ in $R$ and suppose that $\tilde{\iota}_M:M\to \End_R(\tilde{\sF})$ is given by power series $$\tilde{g}_m(x)=\alpha_1(m)\cdot x+\sum_{i\geq 2}\tilde{d}_{m,i}x^i$$
(as has been remarked earlier, note that $m\mapsto\alpha_1(m)$ is a homomorphism of monoids $M\to R$). 
Now define $f:\Z[M][c_{i,j},d_{m,i}]\to R$ as follows: $m\mapsto \alpha_1(m)\in R$ and extend it linearly to  $\Z[M]$. Next map  the variables 
$f(c_{i,j})=\tilde{c}_{i,j}\in R$, and the variables $f(d_{m,i})=\tilde{d}_{m,i}\in R$. Then clearly $g_m(x)\mapsto \tilde{g}_m(x)\in \pow R x$.
Since $\tilde{\sF}$ is an $M$-formal group, its coefficients satisfy all the relevant relations which make it into an $M$ formal group over $R$ (in other words this homomorphism maps the ideal $I$ constructed above to zero), hence this map $f:\Z[M][c_{i,j},d_{m,i}]\to R$ factors through the quotient ring $L_M\to R$ (again denoted by $f$). Then $$\tilde{\sF}(x,y)=x+y+\sum_{i,j} f(\tilde{c}_{i,j})x^iy^j$$ gives rise to $\tilde{\sF}$ through this homomorphism of rings.

Now suppose $f:L_M\to R$ is any homomorphism of rings. Then by applying $f$ to coefficients of the power series $\sF(x,y)$,  and to the coefficients of the endomorphisms $[m]$ (for $m\in M$) of $\sF$, which give the universal $M$-formal group law  $\sF$ over $L_M$, one obtains an $M$-formal group law over $R$.

This proves the assertion.
\ep

\subsection{Functorial properties}
Let me record here the following elementary functorial property which will be used later on.

\bpro\label{prop:functoriality-in-M}
Let $M'\to M$ be a morphism of monoids. Let $\sF_M$ (resp. $\sF_{M'}$) be the universal $M$-formal group (resp. $M'$-formal group) over $L_M$ (resp. $L_{M'}$). Then one has a homomorphism of rings $L_{M'}\to L_M$ such that the $M'$-formal group  $\sF_{M}$ (given by $M'\to M$) arises from the universal $M'$-formal group $\sF_{M'}$ by applying the morphism $L_{M'} \to L_M$ to the coefficients of $\sF_{M'}$. To put it differently, the $M'$-formal group $\sF_M$ is the pull-back of the universal $M'$-formal group $\sF_{M'}$ over $\spec(L_{M'})$ by the morphism $\spec(L_{M})\to\spec(L_{M'})$.
\epro
\bp 
The assertion is immediate from the universal property of $(L_{M'},\sF_{M'})$. The composite $M'\to M\to \End(\sF_M)$ provides an $M'$-formal group over $L_M$.  By the universal property of $(\sF_{M'},L_{M'})$ one has a ring homomorphism $L_{M'}\to L_M$ corresponding to the $M'$-formal group $\sF_{M}$ over $L_M$. This proves the assertion.
\ep
The following is an immediate, but important consequence of the functoriality given by Proposition \ref{prop:functoriality-in-M}:
\bcor\label{cor:isom-univ} 
Let $M$ be a monoid. Let $\sigma:M\to M$ be an isomorphism of monoids. Then $\sigma$ induces an isomorphism $\sigma_M$ of the universal pair $(L_M,\sF_{M})$.
\ecor

\subsection{Existence of a universal $\O^\tr$-formal group and its consequences}
 The following theorem is a consequence of Theorem~\ref{th:universal-M-formal-group}  and Theorem~\ref{th:strict-formal-groups-for fields} as applied to the monoid $\okm$ for a $p$-adic field $K$.

\bthm\label{th:universal-Otr-formal-group}
Let $K$ be  a $p$-adic field and let  $M=\okm$ be the monoid of non-zero elements of the ring of integers $\O_K$.
\benumlab
\item Then there is a universal $\okm$-formal group law over the ring $L_{\okm}$.
\item Let $\sF_1$ be any $\okm$-formal group over a ring $R$ (for example a Lubin-Tate formal group). Then there is a homomorphism of rings $L_{\okm}\to R$ such that $\sF_1$ arises from the universal $\okm$-formal group $\sF$ over $L_{\okm}$.
\eenum
\ethm

\brem\label{re:comparison}\
\benumlab
\item Let me emphasize  that the construction of the ring $L_{\okm}$ and the universal formal group law over $L_{\okm}$  makes no explicit or implicit use of the ring structure of $K$. Indeed, the construction of the monoid ring $\Z[\okm]$ of $\okm$  makes no reference to the additive structure of $K$ (or $\O_K$) and so the construction of Theorem~\ref{th:universal-Otr-formal-group} is independent of the additive structure of $\O_K$ and hence its construction is independent of the ring structures on $\okm\cup\{0\}$.
\item Thus the universal $\okm$-formal group law over $L_{\okm}$ should be viewed as encoding the variation of additive structures on the set $\okm\cup\{ 0\}$ keeping its multiplicative structure fixed and is fully compatible with the ideas set out in \topics, \iut.
\item  
The constructions carried out here should be viewed as being global and more importantly an integral one in the sense that the universal formal group constructed here binds together all $p$-adic field $K$ with multiplicative monoid isomorphic to $\okm$ (the only fixed datum in the construction is the monoid $\O^\tr\isom \okm$). 
\item The universal formal group law provides a (global) access to the additive structure. 
\item In contrast Mochizuki  uses the $p$-adic logarithm to access the additive structure but this necessitates working over  $\Q$ (as the $p$-adic logarithm is defined over $\Q$ and not over $\Z$). In particular the construction is not integral.
\item  On the other hand I ignore the topology of the monoid $\okm$ completely but presumably more refined versions of my constructions could over come this defect.  
\item Nevertheless let me point out that there  is a close resemblance between the two approaches: in Mochizuki's reconstruction theorem for $p$-adic fields  the ring structure of a $p$-adic field emerges from the endomorphism ring of the Lubin-Tate character; and here the ring structure emerges from the endomorphism ring of the formal group.
\eenum
\erem

Let us note one further consequence which is relevant to us: 
\bcor\label{cor:universal-additive-expression} Let $K$ be a $p$-adic field. For any $\alpha\in \okm$, there exists a power series $$[\alpha](x)\in \pow {L_{\okm}} x $$  representing the endomorphism
$$[\alpha]\in \End_{L_{\okm}}(\sF)$$ of the  universal $\okm$-formal group. This power series expresses $\alpha$ independently of the ring structures on $\O_K^\cdas=\okm\cup\{ 0\}$.
\ecor

Also important is the following corollary:

\bcor
Let $K$ be a $p$-adic field. Then any $\O_K$-scheme is a scheme over the scheme $\Spec(L_{\O^\tr})$. In particular schemes $X_1,X_2$, defined over any pair of anabelomorphic fields $K_1,K_2$ (respectively), are schemes over the common base scheme $\Spec(L_{\O^\tr})$ given by Theorem~\ref{th:universal-Otr-formal-group} for the common monoid $\O^\tr\isom\O^\tr_{K_1}\isom \O^\tr_{K_2}$.
\ecor

\bp 
Mochizuki's reconstruction theory provides from $G$ the multiplicative monoid $\O^\tr$ and an isomorphism of monoids $\O^\tr\isom \okm$. This isomorphism together with a choice of a strict $\okm$-formal group over $\O_K$ provides a morphism of rings $L_{\O^\tr}\to\O_K$ equivalently, the arrow $\Spec(\O_K)\to \Spec(L_{\O^\tr})$. So any $\O_K$ scheme can be viewed as a scheme over $\Spec(L_{\O^\tr})$. The remaining part of the assertion follows from this.
\ep

\bpro
Let $M$ be a monoid and suppose that for some ring $R$, there exists an $M$-formal group $\sF$ over $R$ such that $\sF$ is not $R$-isomorphic to the additive formal group $\gah R$. Then the universal $M$-formal group over $\sF_M$ is not $L_M$-isomorphic to $\gah {L_M}$. In particular if $M=\ot$ then $\sF_{\ot}$ is not $L_{\ot}$-isomorphic to $\gah {L_{\ot}}$.
\epro
\bp 
Suppose $\sF_M\isom \gah {L_M}$ over $L_M$. By the universal property of $\sF_M$, there exists a homomorphism of rings $L_M\to R$ such that the $M$-formal group $\sF$ is obtained by applying this homomorphism to the coefficients of $\sF_M$. As $\sF_M\isom \gah {L_M}$ one sees that $\sF\isom \gah R$. But by hypothesis the $M$-formal group $\sF$ is not $R$-isomorphic to $\gah R$. This is a contradiction. 

Now suppose $M=\ot$ and let $K$ be a $p$-adic field such that $\O_K^\tr \isom \ot$. Then by Proposition~\ref{le:lubin-tate-strict} there exists an $\ot$-formal group over $\O_K$ which is not $\O_K$-isomorphic to the additive formal group $\gah {\O_K}$. Thus one sees that $\sF_{\ot}$ is not $L_{\ot}$-isomorphic to $\gah {L_{\ot}}$.
\ep

Let $R$ be a ring and let $\sF$ be a formal group over $R$. I say that $\sF$ is the \emph{trivial formal group over $R$} if $\sF\isom \gah R$ is an isomorphism of formal groups over $R$. If no such isomorphism exists then I say that $\sF$ is a non-trivial formal group over $R$. Preceding result says that the universal $\ot$-formal group $\sF_{\ot}$ is a non-trivial formal group over $L_{\ot}$.
 
Now let me prove the following  converse to Theorem~\ref{th:strict-formal-groups-for fields}(1) which shows that ring structures of interest on $\okm\cup\{ 0\}$ always arise from my constructions:
\bthm\label{th:converse-all-ring-structures} Let $\ot\isom\okm$ be a multiplicative monoid of non-zero elements of the ring of integers  of some $p$-adic field. Let $\ocs=\ot\cup\{ 0\}$ be the multiplicative monoid defined earlier. 
\benumlab
\item Let $R$ be a ring which satisfies the following hypothesis:
\benum[label={\bf(\roman{*})}]
\item Suppose $R$ is a domain of characteristic zero,
\item one has an isomorphism $\ot\isom R^\tr$.
\eenum
 Then the choice of the trivial formal group over $R$ equips $\ocs$ with a ring structure which is isomorphic to the ring $R$.
\item 
Now suppose that $R$ is a ring which satisfies the following hypothesis:
\benum[label={\bf(\roman{*})}]
\item $R$ is a discrete valuation ring of characteristic zero,
\item the residue field of $R$ is finite,
\item one has an isomorphism of multiplicative monoids $\ot \to R^\tr$.
\eenum
Then there is a non-trivial, strict $\ot$-formal group $\sF$ over $R$ such that $\ocs \to \End_R(\sF)\to R$ is an isomorphism of rings.
\item At any rate up to isomorphism, 
\benum[label={\bf(\roman{*})}]
\item Every commutative domain $R$ of characteristic zero such that $\ot\isom R^\tr$ arises from some $\ot$-formal group (and hence from universal pair $(L_{\ot},\sF_{\lot})$) in the manner described by Proposition~\ref{pr:formal-groups-and-ring-struct}.
\item Every discrete valuation ring $R$ with a finite residue field and such that $\ot\isom R^\tr$ arises from some non-trivial $\ot$-formal group (and hence from universal pair $(L_{\ot},\sF_{\lot})$) in the manner described by Proposition~\ref{pr:formal-groups-and-ring-struct}.
\eenum
\eenum
\ethm
\bp Suppose $\ocs=\ot\cup\{0 \}$ is equipped with the structure of a discrete valuation ring of characteristic zero with a finite residue field such that the monoid $\ot\subset\ocs=\ot\cup\{0\}$ is the multiplicative monoid of non-zero elements of this ring. Let us write $R$ for this ring. Let us suppose for the moment that one can find a strict $\ot$-formal group $\sF$ over $R$ such that $\End_R(\sF)\isom R$. Then by the universal property of the pair $(\lot,\sF_{\lot})$, there is a homomorphism of rings $\lot\to R$ such that $\sF$ is the $\ot$-formal group over $R$ obtained by applying $\lot\to R$ to the coefficients of $\sF_{\lot}$. The ring structure on $\ocs$ given by  Proposition~\ref{pr:formal-groups-and-ring-struct} (by virtue of existence of a strict $\ot$-formal group $\sF$ over $R$) is by definition the ring structure of $\End_{R}(\sF)$.  As $\End_R(\sF)\isom R$ is an isomorphism of rings, one has an isomorphism of rings $\ocs\to \End_R(\sF)\to R$.

Thus to complete the proof it suffices to prove that there exists a strict $\ot$-formal group $\sF$ over $R$ with $\End_R(\sF)\isom R$. This is a consequence of our hypothesis that $R$ is a discrete valuation ring of characteristic zero with a finite residue field and of Lubin-Tate theory (see \cite[Section 20.1.19--20.1.21]{hazewinkel-book}).
\ep

\subsection{Existence of Frobenius morphism}
Under natural hypothesis on the monoid $M$, the functoriality of the universal $M$-formal groups provides a Frobenius morphism which is the $p$-th power map at the level of the monoids. Mochizuki's  \cite{mochizuki-frobenioid1} and \iut\ suggest and work with a similar Frobenius operation \textit{but}  purely at a monoid theoretic level.  

\begin{defn}
Let $M$ be a monoid. Let $p$ be a prime number. I will say that $M$ is a $p$-perfect monoid (resp. a perfect monoid) if the homomorphism $\phi_p(m)=m^p$ (resp. $\phi_n(m)=m^n$ for all $n\geq 1$) is an isomorphism  (resp. $\phi_n$ is an isomorphism for all $n\geq 1$) of multiplicative monoids.
\end{defn}

\brem 
Perfect and $p$-perfect monoids play a central role in Mochizuki's theory of Frobenioids \cite{mochizuki-frobenioid1} and also in \iut. 
\erem

\blem 
Let $M$ be a monoid. Then the monoid $$M^{p-pf}=\varprojlim_{n} (M,\phi_{p^n})$$
is $p$-perfect and the monoid $$M^{pf}=\varprojlim_{n} (M,\phi_n)$$ is perfect.
\elem
\bp 
The proof is immediate from the definitions.
\ep

\bthm\label{th:frob-exist}\ 
Let $M$ be a monoid. 
\benumlab 
\item If $M$ is a $p$-perfect monoid, then there exists a natural homomorphism $$\phi_M:L_M\to L_M,$$  called the Frobenius homomorphism,  which is induced by the $p$-th power homomorphism $\phi_p:M\mapright{m\mapsto m^p} M$.
\item If $M$ is a perfect monoid then,  for all $n\geq1$, there exists a natural homomorphisms $\{\phi_n\}_{n\geq1}$,  $$\{\phi_{M,n}:L_M\to L_M\}_{n\geq1},$$ called the Frobenius homomorphisms of $L_M$, which are induced by the $n$-th power homomorphisms $\phi_n:M\mapright{m\mapsto m^n} M$.
\eenum
\ethm
\bp
This is immediate from the functorial properties of the universal $M$-formal group laws Theorem~\ref{th:universal-M-formal-group} and Corollary~\ref{cor:isom-univ}.
\ep

\section{Geometric applications}\label{se:geometric-applications}
\newcommand{\fun}[1]{\Pi_{#1}}
\newcommand{\funt}[1]{\Pi^{{\rm temp}}_{#1}}
In this section I want to consider geometric applications of the theory developed in the preceding sections. The theory of preceding sections can be applied to geometric monoids such as the field of non-zero meromorphic functions on a geometrically connected, smooth, quasi-projective variety over a field. In the context of \topics, \iut, one assumes that this variety is a hyperbolic curve.  
\subsection{Quasi-projective varieties and monoid formal group laws}\label{ss:geometric-case}

Let $K$ be a field and let $G_K$ be its absolute Galois group. While the results of this section  are valid for a broad class of fields, of special interest to me are the following possibilities for $K$: (a) number field, (b) a $p$-adic field. The case $K$ is a perfectoid field is treated in Subsection~\ref{ss:perfectoid-case}.  Let $X/K$ be a geometrically connected, smooth, quasi-projective variety over $K$ ($X$ need not be proper). 

\newcommand{\bkx}{\bar{\mathscr{K}}_X}

Let $\kx:=K(X)$ be the field of rational functions on $X$. Let $\bkx\supset \kx$ be a separable closure of $\kx$. Write $G_X=\gal(\bkx/\kx)$ be the absolute Galois group of $X$. This is obviously a birational invariant of $X$.

\newcommand{\ux}{\mathscr{U}_X}
\newcommand{\bux}{\bar{\mathscr{U}}_X}
Let  $\kx^*$ (resp. $\bkx^*$) be the multiplicative monoid of non-zero elements of $\kx$ (resp. $\bkx$).  Note that one has a natural inclusion of monoids $K^*\into\kxt_X$. Let $\kxs=\kxt_X\cup\{0\}$. Let $$\ux=\Gamma(X,\O_X^*)=H^0(X,\O_X^*)$$ and let $$\bux=\varprojlim_{Y/X} \Gamma(Y,\O_Y^*)$$
where the inverse limit is over all geometrically connected, quasi-projective varieties $Y/L$ with a finite $K$-morphism $Y\to X$ and an inclusion of fields $\ky\subset \bkx$ (since $X$ is quasi-projective but need not be projective, $\ux\supsetneq K^*$ in general). Obviously one has the inclusion of monoids $\ux\subset \kx^*$ and $\bux\subset \bkx^*$.

The formalism of $M$-formal groups developed in preceding sections can also be applied to the monoids $M=\kx^*$ and $M=\bkx^*$ and other monoids obtained from $X$ and one gets the following:

\bthm\label{th:universal-oxt-formal-group}
Let $K$ be a field of characteristic zero and let $X/K$ be a geometrically connected, smooth, quasi-projective variety over $K$. Let $p$ be a prime number. Let $K^*\into \kx^*$ be the inclusion of monoids as above. Let $$M_X\in\left\{ \kx^*,\bkx^*, \bkx^{*p-pf},\bkx^{*pf},\ux,\bux,\bux^{p-pf},\bux^{pf}\right\}.$$ 
\benumlab
\item\label{th:universal-oxt-formal-group-3} There exists a universal $M_X$-formal group $\sF_{M_X}$ over a ring $L_{M_X}$. 
\item If $M_X\in\left\{\bkx^{*p-pf}, \bux^{p-pf}\right\}$, then $(L_{M_X},\sF_{M_X})$ is equipped with a Frobenius endomorphism $\phi_p$ given by Theorem~\ref{th:frob-exist} and if $M_X\in\left\{\bkx^{*pf}, \bux^{pf}\right\}$ then $(L_{M_X},\sF_{M_X})$ is equipped with a family of Frobenius endomorphisms $\{\phi_n \}_{n\geq1}$ given by Theorem~\ref{th:frob-exist}.
\item\label{th:universal-oxt-formal-group-1} As $K$ is of characteristic zero,  strict $\kx^*$-formal groups exist and hence there exists a homomorphism of rings $L_{\kx^*}\to \kx$ whose kernel is a prime ideal.
\item\label{th:universal-oxt-formal-group-2} Any strict $\kx^*$-formal group over a ring $R$ equips $\kx^\circledast$ with a ring structure (given by Proposition~\ref{pr:formal-groups-and-ring-struct}).
\item\label{th:universal-oxt-formal-group-4} The inclusion of monoids $K^*\into\kx^*\into \bkx^*$ provides a ring homomorphisms $$L_{K^*}\to L_{\kx^*}\to L_{\bkx^*}$$ given by Theorem~\ref{prop:functoriality-in-M} making $L_{\kx^*}$ and $L_{\bkx^*}$ into $L_{K^*}$-algebras.
\item\label{th:universal-oxt-formal-group-6} For $M_X\in(\bkx^*,\bkx^{*p-pf},\bkx^{*pf},\bux^{p-pf},\bux^{pf})$ one has an action $$\gal(\bkx/\kx)\act M_X$$ and hence there is a natural action of  $$\gal(\bkx/\kx)\act (L_{M_X},\sF_{M_X})$$  on the pair $(L_{M_X},\sF_{M_X})$.
\eenum
\ethm
\bp 
The assertions {\bf(1)}, {\bf(2)} and  {\bf(6)} are clear from Theorem~\ref{th:exist-M-formal} and Theorem~\ref{th:frob-exist}. The proof of {\bf(4)} is clear from {\bf(1)}, {\bf(2)} and Proposition~\ref{pr:formal-groups-and-ring-struct}. The proof of {\bf(5)} is clear from {\bf(1)}, {\bf(2)} and Theorem~\ref{prop:functoriality-in-M}.

Let me now prove the only remaining assertion {\bf(3)}.
To see that this provides the variation of ring structures on $\kx^*$, I have to show that strict $\kx^*$-formal groups exist. This is done as follows.

For a  $F$ be a field of characteristic zero, let $\gah F$ denote the formal additive group over $F$. Then one checks from Definition~\ref{de:monoid-frml-grp} that $\gah F$ is an $F^*$-formal group over $F$ which is also a strict $F^*$-formal group over $F$ (as $F$ is of characteristic zero). 

Now let $X/K$ be as in the statement of the theorem and take $F=\kx$. Then $F$ is a field of characteristic zero by the hypotheses of the theorem, so one sees that $\gah {\kx}$ is a strict $\kx^*$-formal group over $F=\kx$. Hence strict $\kx^*$-formal groups exist over $F=\kx$. In particular, by Theorem~\ref{th:universal-M-formal-group}, there exists a homomorphism $L_{\kx^*}\to \kx$ whose kernel is a prime ideal since the later is a field (and hence an integral domain).

Another (equivalent way) of constructing $F^*$-formal group is of the following. Fix a line bundle on $X$ and consider its total space as a $\mathbb{G}_a$-bundle over $X$. Let $\sF$ be the formal completion of this $\mathbb{G}_a$-bundle along its zero section. Then $\sF$ is a formal group over $X$. Consider the pull-back of this formal group over $X$ to the generic point of $X$. Then one gets a formal group over the field $F=\kx$. Since $F\supset \Q$ (as $F$ is of characteristic zero), one can trivialize this formal group  to a formal group of the  above form (i.e.  this formal group is $F$-isomorphic to the formal group $\gah F$).

This completes the proof of {\bf(3)} and hence of the theorem.

\ep
\brem
The meaning of Theorem~\ref{th:exist-M-formal}({\bf1, 3}) is that there exists a universal additive law which binds together all the possible additive laws on $\kx^*$ and Theorem~\ref{th:universal-oxt-formal-group}\ref{th:universal-oxt-formal-group-4} establishes the compatibility of addition law on $\kx$  with the addition law on $K$. The assertion Theorem~\ref{th:universal-oxt-formal-group}\ref{th:universal-oxt-formal-group-6} should be compared with the corresponding result in perfectoid geometry proved in Theorem~\ref{th:perfectoid-case}.
\erem

\newcommand{\monopair}[2]{#1\curvearrowright #2}
\newcommand{\mopa}[2]{\monopair{#1}{#2}}
\brem
Let me point out that  the formalism developed here extends naturally to the formalism of  monoid pairs considered by Mochizuki \cite{mochizuki-iut1,mochizuki-iut2}. Let $G$ be a group and $M$ a commutative monoid with an action of $G$. Such a pair will be written as $\mopa G M$ and will be called a monoid pair. Typical example of these pairs which arise in \iut\   (especially \cite{mochizuki-iut2}) are $$\mopa G \bot,\ \mopa G \bos,\ \mopa G \bomu\text{ and }\mopa G {\bK}^*$$ and the geometric monoid pairs $$\mopa {G_X} {\bkx^*},\  \mopa {G_X} {\bkx^{*p-pf}} \text{ and } \mopa {G_X} {\bkx^{*pf}}.$$
The formalism of this paper can be applied to these monoid pairs to obtain universal addition laws for any of these monoid pairs $\mopa G M$.
\erem

\subsection{Perfectoid rings, fields and Huber pairs and monoid formal group laws}\label{ss:perfectoid-case}
One of the important observations of this paper is that theory of perfectoid fields, perfectoid rings and perfectoid Huber pairs is related to the formalism monoid formal group laws. In particular, one may think of the theory of perfectoid fields and rings as a variation of additive structures keeping the multiplicative structure fixed.  

This demonstrates the close relationship between Mochizuki's idea of variation of ring structures keeping multiplicative structures fixed and the ideas in perfectoid geometry of \cite{scholze12-perfectoid-ihes}, \cite{kedlaya-liu15} (the relationship with \cite{fargues-fontaine} is detailed in \cite{joshi-teich-def} and will not be discussed here). This recognition played an important role in the theory of Arithmetic Teichmuller Spaces documented extensively in \cite{joshi-teich,joshi-teich-estimates,joshi-teich-def,joshi-teich-rosetta}. In this section, I record these observations. 

\newcommand{\smcrc}{\circ}\newcommand{\flatp}{\flat\scriptscriptstyle{+}}Let $K$ be a perfectoid field \cite[Definition 3.1]{scholze12-perfectoid-ihes} and let $R$ be a perfectoid $K$-algebra (\cite[Definition 5.1]{scholze12-perfectoid-ihes}). For examples of perfectoid $K$-algebras see \cite[Proposition 5.20]{scholze12-perfectoid-ihes}. Let $K^\flat$ be the tilt of $K$ (\cite[Lemma 3.4]{scholze12-perfectoid-ihes}) and let $R^\flat$ be the tilt of $R$ (\cite[Proposition 5.17]{scholze12-perfectoid-ihes}). Then  $K^\flat$ is a perfectoid field of characteristic $p>0$ and $R^{\flat}\supset K^{\flat}$ and $R^\flat$ is a perfectoid $K^{\flat}$-algebra of characteristic $p>0$. 

Let $R^\smcrc$ (resp. $K^\smcrc$) be the subring of power bounded elements of $R$ (resp. $K$). Note that $K^\smcrc=\O_K$ is the valuation subring of $K$.   

In fact, it is more convenient to work with the more general framework \cite{scholze12-perfectoid-ihes}, \cite{kedlaya-liu15}, \cite{scholze-weinstein-book} of Huber pairs $(R,R^+)$ with $R$ a perfectoid ring and $R^+\subseteq R^\smcrc\subseteq R$ an open, integrally closed subring. By \cite[Lemma 6.2.2]{scholze-weinstein-book}, the ring $R^\flat$ is of characteristic $p>0$. The pair $(R,R^+)$ has tilt $(R^\flat,R^{\flatp})$ and the pair $(R,R^\smcrc)\supseteq (R,R^+)$ has tilt $(R^\flat,R^{\flat\smcrc})$ (resp. $(K^\flat,K^{\flat\smcrc})=(K^\flat,\O_K^{\flat})=(K^\flat,\O_{K^{\flat}})$). Let $W(R^\flat)$ (resp. $W(R^{\flatp})$) be the ring of Witt vectors of $R^\flat$ (resp. $R^{\flatp}$). 

\bthm\label{th:perfectoid-case} 
Let $(R,R^+)$ be a perfectoid Huber pair (e.g. take $(R,R^+)=(K,K^\smcrc)$ for any perfectoid field $K$). Let $(R^\flat,R^{\flatp})$ be the tilt of $(R,R^+)$. Let $M_{R^\flat}=R^{\flat}$ (resp. $M_{R^{\flatp}}=R^{\flatp}$) be the multiplicative monoid of the ring $R^{\flat}$ (resp. $R^{\flatp}$). Then 
\benumlab
\item The monoid $M_{R^\flat}=R^{\flat}$ (resp. $M_{R^{\flatp}}=R^{\flatp}$) is $p$-perfect.
\item There is a universal $M_{R^{\flatp}}=R^{\flatp}$-formal group law $\sF_{R^{\flatp}}$ over the ring $L_{R^{\flatp}}$ and a Frobenius endomorphism $\phi$ of the pair $(L_{R^{\flatp}}, \sF_{R^{\flatp}})$ which arises from the $p^{th}$-power homomorphism of monoid $R^{\flatp}$.
\item Taking $(R,R^+)=(K,K^\smcrc)$, one obtains a universal monoid $K^{\smcrc\flat}$-formal group law $\sF_{K^{\smcrc\flat}}$ over the ring $L_{K^{\smcrc\flat}}$ and a Frobenius endomorphism $\phi$, which coincides with the $p$-power map of the multiplicative monoid $K^{\smcrc\flat}$, of the pair $(L_{K^{\smcrc\flat}}, \sF_{K^{\smcrc\flat}})$.
\item The homomorphism of monoids $K^{\smcrc\flat}\to R^{\flatp}$ provides by functoriality of the universal $K^{\smcrc\flat}$-formal group law (resp. the universal $R^{\flatp}$-formal group law) a homomorphism of rings $L_{K^{\smcrc\flat}}\to L_{R^{\flatp}}$ making $L_{R^{\flatp}}$ into an $L_{K^{\smcrc\flat}}$-algebra.
\item  The tautological $R^{\flatp}$-formal group law $\sF(x,y)=x+y\in W(R^{\flatp})\mybrackets{x,y}$ over $W(R^{\flatp})$ provides, by the universal property of $(L_{R^{\flatp}}, \sF_{R^{\flatp}})$, a homomorphism of rings $L_{R^{\flatp}}\to W(R^{\flatp})$. This homomorphism is compatible with the respective Frobenius homomorphisms of these rings.
\item If $(S,S^+)$ is any untilt of $(R^\flat,R^{\flatp})$ (i.e. one has an isomorphism of the respective tilts $(S^\flat,S^{\flatp})\isom(R^\flat,R^{\flatp})$) then there is a natural homomorphism $L_{R^{\flatp}}\to S^+$  given by the composition 
$$\begin{tikzcd} 
	L_{R^{\flatp}}\ar[r,"\isom"]\ar[rrr,bend left] & L_{S^{\flatp}} \ar[r] & W(S^{\flatp})\ar[r,"\Theta"] & S^{+}.
\end{tikzcd}$$
\eenum
\ethm
\bp 
The assertion that the multiplicative monoid $R^\flat$ (resp.  $R^{\flatp}$) is perfect, is immediate from the fact  
that one has an isomorphism of multiplicative monoids (\cite[Proposition 5.17]{scholze12-perfectoid-ihes})
$$R^{\flat}=\varprojlim_{r\mapsto r^p} R.$$
and (\cite[Lemma 6.2.5]{scholze-weinstein-book})
$$R^{\flatp}=\varprojlim_{r\mapsto r^p} R^{\flatp}$$
respectively.
This proves {\bf(1)}. The assertions {\bf(2)}--{\bf(4)} are clear from earlier results of this paper. To prove {\bf(5)}, it will be enough to prove that the tautological formal group law of $W(R^{\flatp})$ is an $R^{\flatp}$-formal group law. This is seen using the fact that one has a Teichmuller embedding of multiplicative monoids $$R^{\flatp}\mapright{r\mapsto [r]}  W(R^{\flatp}).$$ Now observe that $$x\mapsto [r]\cdot x\in W(R^{\flatp})\mybrackets{x} $$ is an endomorphism of the tautological formal group law $$\sF(x,y)=x+y$$ over $W(R^{\flatp})$. Hence this tautological formal group law is an $R^{\flatp}$-formal group law over $W(R^{\flatp})$ and hence {\bf(5)} follows from the universal property given by Theorem~\ref{th:universal-M-formal-group}.

To prove {\bf(6)}, observe that the first arrow is induced by the isomorphism of the monoids 
$$R^{\flatp}\isom S^{\flatp}$$ which is a consequence of the fact that $(S,S^+)$ has tilt $(S^\flat,S^{\flatp})=(R^\flat,R^{\flatp})$. The second arrow is given by {\bf(5)} for the pair $(S,S^+)$ and the last arrow $\Theta$ is the standard construction of \cite[Lemma 6.2.8]{scholze-weinstein-book}.
\ep

The following is an immediate corollary:
\bcor\label{cor:untilt-binding}
All untilts of a perfectoid Huber pair $(R^\flat,R^{\flatp})$,  of characteristic $p>0$,  arise from the universal formal group law $L_{R^{\flatp}}$ for the monoid $R^{\flatp}$. In particular, if $R^\flat$ is a perfectoid ring of characteristic $p>0$, then every untilt of $R^\flat$ arises from the universal monoid formal group law $(L_{R^\flat},\sF_{R^\flat})$.
\ecor
\bp 
Let $(S,S^+)$ be an untilt of $(R^\flat,R^{\flatp})$. The ring $S^+$ arises in the manner prescribed by Theorem~\ref{th:perfectoid-case}{\bf(6)}. By \cite[Definition 6.1.1]{scholze-weinstein-book}, $S$ is uniform i.e. $S^\smcrc\subset S$ is an open, bounded subring (i.e. $S^\smcrc$ is a ring of definition of $S$ \cite[Lemma 2.2.4]{scholze-weinstein-book}). As $S^+\subset S$ is open by definition and $S^+\subset S^\smcrc$ so $S^+$ is also open and bounded subring of $S$. Hence $S^+\subset S$ is a ring of definition of $S$. 

Again, by \cite[Definition 6.1.1]{scholze-weinstein-book}, $S$ is a Tate ring, so that $S$ contains a topologically nilpotent unit, say $\omega$.  Thus one can invoke \cite[Proposition 2.2.6(2)]{scholze-weinstein-book} to deduce that one has $S=S^+[\omega^{-n}]$ for some $n>0$. 

Thus to determine the pair $(S,S^+)$, it suffices to obtain $S^+$  and this is clear from Theorem~\ref{th:perfectoid-case}{\bf(6)}.  Hence one sees that all untilts of $(R^\flat,R^{\flatp})$ are bound together by the universal $R^{\flatp}$-formal group law  given by Theorem~\ref{th:perfectoid-case}{\bf(2)} for the multiplicative monoid $R^{\flatp}$.
\ep

\brem\  
\benumlab
\item The assertion Theorem~\ref{th:perfectoid-case}{\bf(6)} says that all untilts of $(R^\flat,R^{\flatp})$ are bound together by the universal formal group law for the monoid given in Theorem~\ref{th:perfectoid-case}{\bf(2)}.
\item
In the contexts of arbitrary quasi-projective varieties over $p$-adic or even perfectoid fields, passage to perfectoidification and tilting is not always available  for arbitrary varieties (see \cite[Page 11]{scholze18-icm1}). However, the formalism and results of  Section \ref{ss:geometric-case}, are always available.
\eenum
\erem

\subsection{Global arithmetic analog of Fontaine's $p$-adic ring $\ainf$}\label{ss:global-ainf}
Let $\cpt$ be the tilt of $\C_p$ (\cite{scholze12-perfectoid-ihes}). In the perfectoid setting, and modern $p$-adic Hodge Theory (\cite{fargues-fontaine}, \cite{scholze12-perfectoid-ihes}, \cite{kedlaya-liu15} and others) the ring $\ainf=W(\O_{\cpt})$ introduced in \cite{fontaine94a} plays a fundamental role. I want to remark that the results of Sections \ref{ss:geometric-case}, \ref{ss:perfectoid-case} especially Theorem~\ref{th:universal-oxt-formal-group}, Theorem~\ref{th:perfectoid-case} and Corollary \ref{cor:untilt-binding} suggest a definition of a global i.e. a number field analog $$\ainfl$$ of Fontaine's ring (\cite{fontaine94a}) $$\ainf=W(\cpt)$$ in the $p$-adic setting and a definition of untilts of a number field. Let me provide a brief sketch of the construction of this global ring (detailed study of the properties of $\ainfl$  and its relation with Mochzuki's global arithmetic Kodaira-Spencer classes (\cite{mochizuki-HA})  is taken up in \cite{joshi-teich-ks-theory} which is under preparation).

For a ring $R$, let $\bbW(R)$ be the ring of big Witt-vectors of $R$ (\cite{hazewinkel-book}); let $L$ be a number field, let $\bL$ be an algebraic closure, let $\O_L$ (resp. $\O_{\bL}$) be the ring of integers of $L$ (resp. $\bL$). Let \be\label{eq:perfection-M} M=\O_{\bL}^{pf}=\invlim_{n\geq1}(\O_{\bL}, x\mapsto x^n)\ee be the perfection of the multiplicative monoid of the ring of integers $\O_{\bL}$ of $\bL$. This is the notion of perfection of monoids introduced in \cite{mochizuki-frobenioid1} and hence this means $M$ is a perfect multiplicative monoid in the sense of \cite{mochizuki-frobenioid1}. 

Note that  one has a natural action $G_L\act M$ of the absolute Galois group $G_L$ on $M$, as well as  surjective  Frobenius morphisms $x\mapsto x^n$ for all $n\geq1$. An element of $M$ can described as a compatible system of $n^{th}$-roots (all $n\geq 1$) of some $x\in\O_{\bL}$. 

Let $\Z[M]$ be the monoid ring of $M$.  Then Theorem~\ref{th:universal-oxt-formal-group}, Theorem~\ref{th:perfectoid-case} and Corollary \ref{cor:untilt-binding} suggest that \be\label{eq:aninl-def}\ainfl=\bbWl=\bbW(\Z[M])\tensor_{\Z}\O_L=\bbW(\Z[\O_{\bL}^{pf}])\tensor_{\Z}\O_L\ee equipped with 
\benumlab
\item a surjective homomorphism of rings $$\bbWl(\Z[\O_{\bL}^{pf}])\to \O_{\bL},$$ 
(so the kernel is a prime ideal $P=\ker(\bbWl(\Z[\O_{\bL}^{pf}])\to \O_{\bL})$),
\item the given structure \eqref{eq:aninl-def} of $\O_L$-algebra on $\bbWl(\Z[\O_{\bL}^{pf}])$,
\item the Teichmuller-representative homomorphism of a perfect multiplicative monoid 
$$M=\O_{\bL}^{pf}\to \bbWl(\Z[M])=\bbWl(\Z[\O_{\bL}^{pf}])$$
such that the composite homomorphism of multiplicative monoids $$M=\O_{\bL}^{pf}\to \bbWl(\Z[\O_{\bL}^{pf}])\to \O_{\bL}$$ induces an isomorphism on the passage to perfections of monoids;
\item   Frobenius endomorphisms  induced from the Frobenius endormorphisms $x\mapsto x^n$ (for each $n\geq 1$) of $\bbW(\Z[\O_{\bL}^{pf}])$;
\item the canonical action of the absolute Galois group $G_L$ (through its action on $M$),
\eenum
is the global analog of $\ainf$. The set of ``untilts of $\bL$'' should be a Galois, Frobenius stable subset of $\Spec(\bbWl(\Z[\O_{\bL}^{pf}]))$ modulo the action of $\O_L^*$ (by \eqref{eq:aninl-def}). This ring is studied in \cite{joshi-teich-ks-theory}.

The following assertion is now clear:
\bthm
Let $L$ be a number field, $\bL$ and an algebraic closure of $L$, $M$ be the perfection of the multiplicative monoid $\O_{\bL}$ defined in eq.~\eqref{eq:perfection-M} and  let the other notations of this subsection be in force. Then there is a natural homomorphism $$L_M\to\ainfl,$$ from the universal ring $L_M$ provided by the universal formal group law of Theorem~\ref{th:universal-oxt-formal-group}  for the monoid $M$, to the ring  $\ainfl$ defined by eq.~\eqref{eq:aninl-def}.
\ethm
\bp 
The proof is clear from  the proofs of Theorem~\ref{th:universal-oxt-formal-group} and Theorem~\ref{th:perfectoid-case}.
\ep
\bibliographystyle{plainnat}

\end{document}